\documentclass[12pt,a4paper,oneside,notitlepage]{amsart}
\usepackage{amssymb}
\usepackage[arrow,matrix,curve]{xy}
\def\xyma{\xymatrix@M.7em}
\usepackage[centertags]{amsmath}
\usepackage[T1]{fontenc}
\usepackage[latin1]{inputenc}
\usepackage{amsfonts}
\usepackage{amscd}
\usepackage{euscript}
\usepackage{amsmath}

\numberwithin{equation}{section}
\newtheorem{theorem}{Theorem}

\newtheorem{prop}[theorem]{Proposition}
\newtheorem{lemma}[theorem]{Lemma}

\usepackage{amsfonts,amssymb}

\newcommand{\mono}{\hookrightarrow}

\newcommand{\nc}{\newcommand}

\nc{\Ker}[1]{\mbox{Ker$(#1)$}}

\nc{\impl}{\Rightarrow}

\makeindex
\begin{document}
\newcommand{\N}{\noindent}

\title{A ``working mathematician's'' definition of  semi-abelian categories}


\author{Manfred Hartl \& Bruno Loiseau}
\maketitle

\begin{center}\small Univ Lille Nord de France, F-59000 Lille, France\\
UVHC, LAMAV and FR CNRS 2956, 
F-59313~Valenciennes, France
\end{center}

\begin{abstract}
Semi-abelian and finitely cocomplete homological 
categories are characterized in terms of four resp.\ three simple axioms, in terms of the basic categorical notions introduced in the first few chapters of MacLane's classical book. As an immediate application we show that categories of diagrams in semi-abelian and similar categories are of the same type; in particular, the category of  simplicial or $\Gamma$-objects in a semi-abelian category is semi-abelian.

\end{abstract}
\vspace{5mm}

The notion and theory of abelian categories play a crucial role in homological algebra, algebraic geometry and many other fields for a long time. From the very beginning, however, various attempts were made to establish a convenient generalization which would include the category of groups, as it shares many key features with abelian categories, 
in particular most properties of  exact sequences including all classical diagram lemmas. Nevertheless, a satisfactory solution of this problem was achieved only around 
the early 2000's,  in establishing the notions and theory of homological and semi-abelian categories (see \cite{JMT},  \cite{BB}). These notions are powerful enough to ensure a maximum of desirable properties (in particular, all diagram lemmas of homological algebra), but are also general enough to cover a maximum of interesting examples: all categories of algebraic objects having a group law as part of the structure and admitting a zero object, are semi-abelian. This includes the categories of groups and of Lie algebras, and more generally, the category of algebras over any reduced operad. But also the categories of crossed modules, of $\mathbb{C}^*$-algebras and of compact Hausdorff-spaces are semi-abelian, and so are categories of diagrams with values in a semi-abelian category (see Theorem \ref{expo} below). In particular, this includes  categories of simplicial or $\Gamma$-objects in any of the foregoing categories, which crucially occur in algebraic topology. 
A comprehensive and well written account of the fundamental theory of semi-abelian and similar categories  is given in \cite{BB}. However, the very definition of a semi-abelian category as being  a pointed,  exact and protomodular category with 
finite sums, may appear awkward to a ``working mathematician''. So the aim of this paper is not to contribute to the development of the theory itself, but to remedy this unconvenience by providing a characterization of homological and semi-abelian categories in very basic terms, and thus to render these notions more attractive for a wider public.
This seems to be desirable as not only the authors are convinced that the theory of semi-abelian categories will play a similarly fundamental role in future mathematics as the theory of abelian categories does nowadays. For example, important foundations of semi-abelian homological and homotopical algebra are layed  in the work of Van der Linden and others, see e.g. \cite{Hopf}, \cite{simhomt}, \cite{sat}.

To state our main result right away consider the following axioms about a category $\mathbb C$, where we denote by kernel or cokernel the corresponding injection or projection arrow, resp.\vspace{4mm}

\N{\bf A1.} $\mathbb C$ is pointed and has finite products and coproducts, and equalizers and coequalizers; hence it has all finite limits and colimits;\medskip

\N{\bf A2.} For any split epimorphism $p : X \rightarrow Y$ with section $s : Y \rightarrow X$ and with kernel $k : K[p] \mono X $, the arrow $<k, s>: K\amalg Y\rightarrow X$ is a cokernel;\medskip

\N{\bf A3.} The pullback of a cokernel  is a cokernel;\medskip

\N{\bf A4.} The image of a kernel by a cokernel is a kernel.\vspace{4mm}

Note that when  $\mathbb C$ is the category of groups,  axiom A2 says that a semi-direct product group $G=N \rtimes T$ is generated by $N$ and $T$.

Recall from \cite{BB} that roughly speaking, homological categories are designed to make all the lemmas of homological algebra hold, while semi-abelian categories are homological but also allow for additional constructions such as, notably, semi-direct products and crossed modules, see also \cite{BJ}, \cite{J}. \medskip


\N{\bf Theorem.}
\textit{A category $\mathbb C$ is  finitely cocomplete homological if and only if it satisfies the three axioms A1-A3, and is semi-abelian   iff it satisfies all the four  axioms A1-A4.}\medskip

The proof is achieved at the end of section 2 below.
\medskip

\N{\bf Remarks.}

1) The notion of a kernel or a cokernel can be defined in any pointed category as (co)equalizers; hence axioms A1-A3 are indeed stated in categorical terms; also note that a cokernel always is the cokernel of its kernel if the latter exists (here assured by axiom A1). The notion of an image, however, requires the setting of {\it regular} categories, where  any morphism $f$ can be factored (uniquely up to isomorphism) as a regular epimorphism followed by a monomorphism (a regular epimorphism is a morphism which is a coequalizer). This monomorphism, considered as a subobject, is the image of $f$. And it will be shown below that categories satisfying Axioms A1 and A2 (or even some weakened form of these axioms) and A3 are regular ; so Axiom A4 makes sense when Axioms A1, A2 and A3 hold. 

2) The axioms are stated in terms of cokernels, but at every place in these axioms ``cokernel'' may be replaced by ``regular epimorphism'', which is more common in papers in this field. In fact, it follows from Lemma \ref{Manfred} and Proposition \ref{prodsemdir} below that in a category   satisfying Axiom A1 and  a weaker form of Axiom A2 in which ``cokernel'' is replaced by ``regular epimorphism'',  any regular epimorphism is a cokernel.


Moreover, a homological category is finitely cocomplete if and only if it has finite sums and  coequalizers of (internal) equivalence relations. For a given category, existence of the latter often is easily checked (an internal equivalence relation in many ``concrete'' categories essentially is a congruence), so the main difficulty in proving (finite) cocompleteness consists of the proof of existence of (finite) sums. This characterization of finite cocompletess is explicitly shown in \cite{BB} (though expressed in a slightly different way, because of the exactness condition),  in the stronger case of semi-abelian categories, but  just a few words have to be changed in the proof to get this slightly stronger result. We point this out because interesting examples of finitely cocomplete homological categories which are not semi-abelian do actually exist (e.g. the category of topological groups, or the category of pairs of groups in section 1 below). 

Finally, we advertize that in forthcoming work the authors study special semi-abelian categories called nilpotent
 (categories of nilpotent groups or algebras being the guiding examples, and where abelian categories are just those of nilpotency class 1);
this provides a categorical foundation for quadratic algebra (and of non-linear algebra of higher degree, later on) which was inaugurated in \cite{BHP} and is further developped by Baues, Jibladze, Pirashvili, Muro, and more recently, also by Gaudier, Goichot and the authors,  see for instance the references [1]-[6],  \cite{P2RM} and \cite{QMaps}. This theory   will be presented to a wider public in a forthcoming book on quadratic algebra jointly written by several of the above authors.


\section{Finitely complete protomodular pointed categories}
We start by recalling some elementary facts; carefully explained details about these notions can be found in \cite{BB}.
All the categories considered here are finitely complete. Recall that a {\it strong epimorphism} is a morphism $q:X\rightarrow Y$ such that if it factors through a monomorphism $m:Z\mono Y$, then this monomorphism is an isomorphism. In a finitely complete category, a strong  epimorphism indeed is an epimorphism, and any regular epimorphism is strong. A morphism which is both a monomorphism and a strong epimorphism is an isomorphism.
Also recall that a  family $(f_i:X_i\rightarrow Y)$ of morphisms with same codomain is {\it epimorphic} if and only if for any pair of morphisms $u$, $v$ :  $Y\rightarrow Z$, such that $uf_i=vf_i$ for all $i$, one has $u=v$. Intuitively, the family ``covers'' $Y$. A {\it strongly epimorphic family} is a family $(f_i:X_i\rightarrow Y)$ of morphisms with same codomain such that if all $f_i$'s factorize through the same monomorphism $m:Z\mono Y$, then this monomorphism is an isomorphism.  For instance, the injections of the vertices of a diagram into its colimit (when it exists) form a strongly epimorphic family. Here again, in a finitely complete category, a strongly epimorphic family is epimorphic.

Finally recall that a {\it pointed category} is a category with a zero object, i.e. an object $0$ which is both initial and final. Then for any two objects $X, Y$ there exists a unique morphism $0_{XY}$ which factors through $0$; obviously one has $0_{YZ}0_{XY}=0_{XZ}$. That's why one usually omits the indices. The {\it kernel} of a morphism then is the equalizer of this morphism and the zero arrow, and its {\it cokernel} is,  when it exists, their coequalizer. The notions of epimorphism, strong epimorphism, regular epimorphism and cokernel do not coincide in general; but in homological categories, the latter three do.\\

We give  two examples of categories in which not any epimorphism is strong; they also are of interest, because they turn out to be complete and cocomplete homological, but not semi-abelian.

The category of topological groups (or more generally of topological $\mathbb T$-algebras, where $\mathbb T$ is a protomodular algebra), is homological, complete and cocomplete (see \cite{BB}). However in this category (which is not exact hence not semi-abelian), any surjective continuous map is an epimorphism, while regular morphisms (quotients) are open surjections. 

The following category is another example, which plays an important role in quadratic algebra.\\
Consider the category whose objects are ordered pairs $(G,A)$ where $G$ is a (necessarily 2-step nilpotent) group and $A$ a subgroup of $G$ such that $G'\subseteq A\subseteq \mbox{Z}(G)$, where $G'$ and $ \mbox{Z}(G)$ respectively denote the derived subgroup and the center of $G$. A morphism $f$ between two such objects $(G, A)$ and $(H, B)$ is a group morphism $f$ between $G$ and $H$ such that $f(A)\subseteq B$. It is obvious that  any morphism which, simply considered as a group morphism, is an epimorphism, also is an epimorphism in this category. But such a morphism is not always regular. Indeed, consider two morphisms $f,g : (G,A)\rightarrow(H,B)$ in this category. It is easy to compute their coequalizer : if $q: B\rightarrow Q$ is their coequalizerin the category of groups, then $(Q, q(B)$ is an object in our category and $q$ obviously is a morphism $(G,B)\rightarrow(Q,q(B))$ which is the coequalizer of $f$ and $g$. So, forgetting about $(G,A)$, if a morphism $q : (H,B)\rightarrow (Q,C)$ is a regular morphism in our category, then not only $q$ is a group epimorphism, but also $q(B)=C$. Note that since any group epimorphism is regular (even is a cokernel), this necessary condition is also a sufficient one. So if $G$ is a 2-step nilpotent group in which $G'\subsetneqq \mbox{Z}(G)$, then  the morphism $1_G : (G,G')\rightarrow (G, \mbox{Z}(G))$ is an epimorphism which is not strong.

Of course, the notion of kernel (and of cokernel) is of no interest in general pointed categories, because it carries not enough information (consider for instance the case of the category of pointed sets). The protomodularity condition, however, which will turn to be equivalent (for a finitely complete pointed category) to condition (A) in the following lemma, is strong enough to make these notions interesting.

 We do not recall the notion of a protomodular category here, but just recall that Proposition 3.1.2. of \cite{BB} states that a pointed category is protomodular if and only if it has pullbacks of split epimorphisms along any map, and if the "split short five lemma" holds, meaning that in the following commutative diagram with $\Ker{q}=i$, $\Ker{p}=k$, $qr=1_W$, $ps=1_Y$, $nq=pm$, $mr=sn$ and $mi=kl$, if the outer vertical arrows $l$ and $n$ are isomorphisms, then so is the central one : 
 
$$\xymatrix{1\ar[r]&K[q]\ar[d]^l\ar[r]^i&Z\ar@<-2pt>[r]_q\ar[d]^m&W\ar@<-2pt>[l]_r\ar[d]^n\ar[r]&1\\
1\ar[r]&K[p]\l\ar[r]^k&X\ar@<-2pt>[r]_p&Y\ar@<-2pt>[l]_s\ar[r]&1}$$

\begin{lemma}
\label{Manfred}

Let $\mathbb C$ be a finitely complete pointed category. Consider the following conditions :

(A) For any split epimorphism $p : X \rightarrow Y$ with section $s : Y \rightarrow X$ and with kernel $k : K[p] \mono X$, the pair $(k, s)$ is a strongly epimorphic family.

(B) Let $f, g : A \rightrightarrows B$ be a parallel pair of epimorphisms with common section $s : B \rightarrow A$. Let $k : K[f] \rightarrow A$ be the kernel of $f$. Then a morphism $q : B \rightarrow C$ is a coequalizer of $f$ and $g$ if and only if it is a cokernel of $gk$. 

In other (but less precise) terms, "the coequalizer of $f$ and $g$ is the cokernel of the restriction of $g$ to $Ker (f)$"

(C) Any morphism whose kernel is $0$ is a monomorphism.

(D) For any  commutative diagram of the following type, if $i= \Ker{q}$ and if  $l$, $k$ and $n$ are monomorphisms, then so is $m$.

$$\xymatrix{%
K[q]\ar[d]^l\ar@{^{(}->}[r]^i&Z\ar[d]^m\ar[r]^q&W\ar[d]^n\\%
K\ar@{^{(}->}[r]^k&X\ar[r]^p&Y%
}$$

(Of course $kl$ can be replaced by any monomorphism $K[q]\mono X$, but we present the diagram in this way to show the analogy with the short five lemma situation where moreover $k$ is a kernel and $q, p$ are (split) regular epimorphisms).

(E) Any regular epimorphism is a cokernel.

Then the following implications hold in $\mathbb C$ : $(A) \impl (B) \impl (C) \Leftrightarrow (D)$ and $(B) \impl (E).$\end{lemma}

\begin{proof} 

{\it (A) implies (B)}: Consider $f, g, s, k$ as in condition (B). Then for a morphism $h : B \rightarrow D$ we have $hf = hg$ if and only if $hgk = 0$. Indeed, $(k,s)$ being a (strong) epimorphic family, $hf=hg$ is equivalent to $hfk = hgk$ and $hfs =hgs$. Since $fs = gs = 1_B$ and $fk=0$, this is equivalent to $hgk=0$. Therefore, a morphism $q$ is a coequalizer of $f$ and $g$ if and only if it is a cokernel for $kg$.

{\it (B) implies (C)}: Let $f : A \rightarrow B$ be a morphism whose kernel is the null object. Consider $p_1, p_2 : P[f]  \rightarrow A$ be the kernel pair of $f$. It is well-known in general that for any morphism $f$, the kernel of $f$ also is the kernel of $p_2$, more precisely that if $k : K[f] \mono A$ is the kernel of $f$, then the only factorization $l : K[f]\rightarrow P[f]$ such that $p_2 l = 0$ and $p_1 l = k$, which exists because $fk = f0 =     0$, is the kernel of $p_2$. Thus here the kernel of $p_2$ is $0 : 0\rightarrow P[f]$. But $p_1$ and $p_2$ are two parallel epimorphisms with common section $\delta$ (the diagonal arrow). Hence by (B), since the identity of $A$ is a cokernel of $0\rightarrow A$, it is also a coequalizer of $p_1$ and $p_2$, which classically characterizes the fact that $f$ is a monomorphism.

{\it (C) is equivalent to (D)}: Consider a diagram as in (D). By (C), we have to show that if $u:U \rightarrow Z$ is such that $mu=0$, then $u=0$. But if $mu=0$, then $pmu=nqu=0$, hence $qu=0$ since $n$ is a monomorphism. So $u$ factors through $K[q]$, say $u=iv$. Then $klv=miv=mu=0$, hence since $kl$ is a monomorphism, $u=0$ as required. Conversely, since in any pointed category any map with source 0 is a monomorphism, the following diagram shows that if $0$ is the kernel of a map $f$ and if (D) holds, then $f$ is a monomorphism.
$$\xymatrix{%
0\ar@{=}[d]\ar@{^{(}->}[r]^0&Z\ar[d]^f\ar[r]^f&W\ar@{=}[d]\\%
0\ar@{^{(}->}[r]^0&W\ar@{=}[r]&W%
}$$

{\it (B) implies (E)}: Consider a morphism $q$ which is the coequalizer of a parallel pair of morphisms $\raisebox{.7ex}{\xymatrix{X{\ar@<2pt>[r]^f\ar@<-2pt>[r]_g}&Y}}$. Then it is well-known that $q$ is also the coequalizer of its kernel pair $(p_1,p_2)$ which is the pullback of $q$ along itself. But this kernel pair admits the diagonal as a common section, hence by (B) its coequalizer also is a cokernel.\end{proof}
\begin{prop}\label{critproto}

The following properties 1) - 3) are equivalent for a finitely complete pointed category $\mathbb C$ :

1) $\mathbb C$ is protomodular.

2) For any split epimorphism $p : X \rightarrow Y$ with section $s : Y \rightarrow X$ and with kernel $k : K[p] \mono X$, the pair $(k, s)$ is a strong epimorphic family (Condition (A) above).

3) a) $\mathbb C$ satisfies the "split short five lemma for strong epis", meaning that in the following commutative diagram with $\Ker{q}=i$, $\Ker{p}=k$, $qr=1_W$, $ps=1_Y$, $nq=pm$, $mr=sn$ and $mi=kl$, if the outer vertical arrows $l$ and $n$ are strong epis (not necessarily isomorphisms), then so is the central one :
$$\xymatrix{1\ar[r]&K[q]\ar[d]^l\ar[r]^i&Z\ar@<-2pt>[r]_q\ar[d]^m&W\ar@<-2pt>[l]_r\ar[d]^n\ar[r]&1\\
1\ar[r]&K[p]\l\ar[r]^k&X\ar@<-2pt>[r]_p&Y\ar@<-2pt>[l]_s\ar[r]&1}$$

and b)  in $\mathbb C$ any morphism whose kernel is the zero object, is a monomorphism (condition (C) above).

\end{prop}

\begin{proof}

{\it1) implies 2)}:

Let $\mathbb C$ be a protomodular category and  $p : X \rightarrow Y$ be a split epimorphism with section $s : Y \rightarrow X$ and with kernel $k : K[p] \mono X$. To prove that the pair $(k, s)$ is a strong epimorphic family, consider a monomorphism $m : Z \rightarrow X$ and two factorizations $f:K[p]\rightarrow Z$ and $t : Y\rightarrow Z$ of $k$ and $s$ through $m$, hence $mf = k$ and $mt=s$. We have to show that $m$ is an isomorphism. But under these conditions $pm : Z \rightarrow Y$ is split epi, with section $t$, since $pmt=ps=1_Y$. Moreover $1_Y( pm)=(pm)$, $mt=s1_Y$ and $mf=1_K k$. We now claim that $f$ is the kernel of $pm$. Indeed, first of all $pmf=0$. Secondly, $f$ obviously is a monomorphism. And if $u:W \rightarrow Z$ is such that $pmu=0$, then since $k$ is the kernel of $p$, $mu$ factors through $k$, say $mu=kv$. But then $mfv=kv=mu$, so since $m$ is a monomorphism, $fv=u$. This gives the (unique) factorization as required.
So we may apply the split short five lemma to the following situation :

$$\xymatrix{1\ar[r]&K[p]\ar[r]^f\ar@{=}[d]&Z\ar@<-2pt>[r]_{pm}\ar[d]^m&Y\ar@<-2pt>[l]_t\ar@{=}[d]\ar[r]&1\\
1\ar[r]&K[p]\ar[r]^k&X\ar@<-2pt>[r]_p&Y\ar@<-2pt>[l]_s\ar[r]&1}$$
to conclude that $m$ is an isomorphism as required.

{\it 2) implies 3a)}: Consider the following commutative diagram where  $\Ker{q}=i$, $\Ker{p}=k$, $qr=1_W$, $ps=1_Y$, $nq=pm$, $mr=sn$ and$mi=kl$, and the outer vertical arrows $l$ and $n$ are strong epimorphisms. We have to show that $m$ is a strong epimorphism as well, i.e. that if $m$ factors through some monomorphism $u:T\mono X$, then this monomorphism is an isomorphism.

$$\xymatrix{1\ar[r]&K[q]\ar[d]^l\ar[r]^i&Z\ar@<-2pt>[r]_q\ar[d]^m&W\ar@<-2pt>[l]_r\ar[d]^n\ar[r]&1\\
1\ar[r]&K[p]\l\ar[r]^k&X\ar@<-2pt>[r]_p&Y\ar@<-2pt>[l]_s\ar[r]&1}$$

Because of condition 2), $(k,s)$ is a strongly epimorphic family. But then, since $l$ and $n$ are strong epimorphisms, $(kl,sn)$ also is a strongly epimorphic family. And if $m$ factors through the monomorphism $u$, then so do $mi$ and $mr$, hence also $kl$ and $sn$ ; and since these form a strongly epimorphic family, $u$ is an isomorphism.

Of course 2) implies 3b), as shown in Lemma 1.

{\it 3) implies 1)}: in the situation of the split short five lemma with the outer vertical arrows being isomorphisms, the central arrow is a mono because of (D) in Lemma 1 (which is equivalent to (C)), and is a strong epimorphism because of the "split short five lemma for strong epis". Hence it is an isomorphism as required.\end{proof}

\section{Homological categories, finitely cocomplete homological categories, semi-abelian categories}

\begin{prop}\label{prodsemdir}A finitely complete pointed category with finite sums is protomodular if and only if it satisfies the following condition

(F) For any split epimorphism $p:X \rightarrow Y$  with section $s$ and kernel $k: K\rightarrow X$, the factorization $(k,s):K \coprod Y\rightarrow X$ is a strong epimorphism. \end{prop}

\begin{proof}Indeed, it is a general fact that in a, say finitely complete category with finite sums, a family $(f_1, f_2)$ of morphisms $X_i \rightarrow X$ is strongly epimorphic if and only if the factorization $(f_1,f_2) : X_1\coprod X_2 \rightarrow X$ is a strong epimorphism.\end{proof}

We now are interested in homological categories, and an particular finitely cocomplete ones. Recall that according to \cite{BB}, a homological category is a pointed, regular and protomodular category. Regularity means that the category is finite complete, that the pullback of any regular morphism along any map is a regular morphism, and that kernel-pairs have a coequalizer (the kernel pair of a morphism is the parallel pair formed by the two projection in the pullback of the morphism along itself)

Notice that in Proposition 5.1.3. of \cite{BB}, where it is proved that every semi-abelian category $\mathbb C$ is finitely cocomplete, they only use the fact that the category is homological with finite sums, and  coequalizers of equivalence relations do exist in $\mathbb C$. Indeed, given a pair of parallel morphisms, they construct an equivalence relation whose coequalizer is the same, if it exists, as the coequalizer of this pair, and they use  exactness to conclude that this coequalizer exists because an equivalence relation is a kernel pair. So the full exactness is not really needed, existence of coequalizers of equivalence relations suffices. This is an interesting fact, because there are examples of categories which are homological and have sums and (easily computable) coequalizers of equivalence relations, but which are not exact (for instance, topological vector spaces, but also modules over a square ringoid (see \cite{P2RM}).

So we have proved :

\begin{prop}\label{homolcoc}A category is homological and finitely cocomplete if and only if it is pointed, finitely complete, satisfies condition (F) in Proposition \ref{prodsemdir}, has finite sums, has coequalizers of equivalence relations, and pullbacks of regular epimorphisms are regular epimorphisms.\end{prop} Note that the two last conditions express slightly more than regularity, but significantly less than exactness.

Finally, we provide a characterization of exactness in the case of a homological category:

\begin{prop}\label{homolex}A homological category is exact if and only if in this category, equivalence relations have a coequalizer and the image of a kernel by a regular epimorphism itself is a kernel.\end{prop}
\begin{proof}

The condition is necessary : existence of coequalizers of equivalence relations in exact categories is obvious. The second condition is an immediate consequence of Propositions 3.2.7. and 3.2.20 of \cite{BB}.

Conversely, suppose that the conditions are satisfied; we show that the category is exact.

  Let $\raisebox{.7ex}{ \xymatrix{R {\ar@<2pt>[r]^{r_0}\ar@<-2pt>[r]_{r_1}}&X\ar@/_1pc/[l]_{\delta}}}$ be an equivalence relation, we must prove that it is a kernel pair relation. By hypothesis, it has a coequalizer $q:X\rightarrow Q$ which itself has a kernel pair (provided with its diagonal $\sigma$) : $\raisebox{.7ex}{ \xymatrix{P{\ar@<2pt>[r]^{p_0}\ar@<-2pt>[r]_{p_1}}&X\ar@/_1pc/[l]_{\sigma}}}$. Since $P$ (or more precisely the pair $(p_0,p_1)$) is the kernel pair of $q$ and $qr_0=qr_1$, there exists a unique morphism $i:R\rightarrow P$ such that $p_0i=r_0$ and $p_1i=r_1$ as shown in this diagram :
  
$$
  
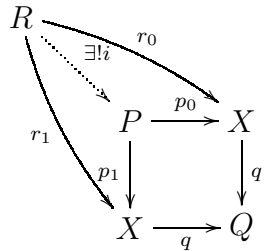
\xymatrix{R \ar@/^/[rrd]^{r_0} \ar@/_/[rdd]_{r_1} \ar@{.>}[rd]^{\exists !i}&&\\&P \ar[r]^{p_0} \ar[d]_{p_1}&X\ar[d]^q\\&X\ar[r]_q&Q}$$or  : 
  
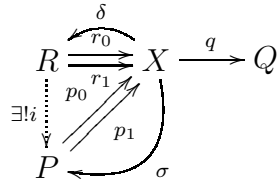
$$\xymatrix{R \ar@{.>}[d]_{\exists !i}\ar@<2pt>[r]^{r_0}\ar@<-2pt>[r]_{r_1}&X\ar@/_1pc/[l]_{\delta}\ar@/^2pc/[ld]^{\sigma}\ar[r]^q&Q\\%
  P\ar@<2pt>[ru]^{p_0}\ar@<-2pt>[ru]_{p_1}&&}$$
  
So to insure that $(r_0,r_1)$ is a kernel pair it suffices to prove that $i$ is  mono and a regular epi, hence an isomorphism.
  
{\it The arrow $i$ is a monomorphism}. Consider two parallel morphisms $a$, $b$ with target $R$. We have to show that if $ia=ib$ then $a=b$. But since $(r_0,r_1)$ is an equivalence relation, it is a monomorphic family, so it suffices to show  that $r_0a=r_0b$ and  $r_1a=r_1b$ which means  $p_0ia=p_0ib$ and  $p_1ia=p_1ib$, which is obvious.

{\it The arrow $i$ is a regular epimorphism}.

Let \xymatrix{\raisebox{1.4ex}{$K$}\ar@{^{(}->}[r]^k&\raisebox{1.4ex}{$R$}} be the kernel of $r_0$ and let $\xymatrix{\raisebox{.7ex}{$K$}\ar@{->>}[r]^p&\raisebox{.7ex}{$L$}\ar@{^{(}->}[r]^l&\raisebox{.7ex}{$X$}} $be the regular epi-monomorphism decomposition of $r_1k =p_1ik$ (in other terms, $L$ is the image by $r_1$ of the subobject $K$  of $R$). Consider the following diagram :

$$
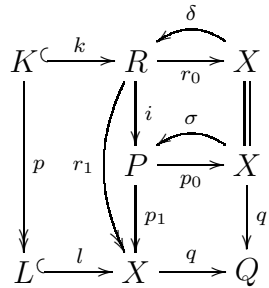
\xymatrix{K\ar@{^{(}->}[r]^k \ar@{->>}[dd]^p&R\ar[r]_{r_0}\ar[d]^i\ar@/_1pc/@{->>}[dd]_{r_1}&X\ar@/_1pc/[l]_{\delta}\ar@{=}[d]\\%
&P\ar[d]^{p_1}\ar[r]_{p_0}&X\ar[d]^q\ar@/_1pc/[l]_{\sigma}\\%
L\ar@{^{(}->}[r]^l &X\ar[r]^q&Q}$$

The only property  to show to prove its commutativity is $\sigma=i\delta$. But since $P$ is a pullback, it suffices to show that $p_j\sigma = p_ji\delta$ for $j=0,1$, which is immediate since one has  $p_ji\delta=r_j\delta=1_X=p_j\delta$.

  Since both $p$ and $r_1$ are regular epimorphisms, $l$ is the image of the kernel $k$ hence is itself a kernel. On the other side, $q$ is the coequalizer of $r_0$ and $r_1$, who have a common section $\delta$. Hence, because of (B) in Lemma \ref{Manfred}, $q$ is also the cokernel of $r_1k$, since $k$ is the kernel of $r_0$. So, $q$ is the cokernel of $lp$, hence of $l$ since $p$ is a (regular) epimorphism.
  
  Since $L$ is the kernel of $q$ and $P$ is a pullback, it is well-known that $l$ also is (isomorphic to) the kernel of $p_0$ ; more precisely, noting that in the following diagram $q0=ql$, the unique factorization $\phi$ of $0$ and $l$ by $p$ is the kernel of $p_0$ : 
  
  $$\xymatrix{&P\ar[r]^{p_0}\ar[d]^{p_1}&X\ar[d]^q\\%
  L\ar@/^4pc/[rru]^0\ar@{.>}[ru]^{\exists !\phi}\ar@{^{(}->}[r]^l &X\ar[r]^q&Q}$$
  Then we may apply the split short five lemma for strong epimorphisms to the following situation to conclude that $i$ is a strong epimorphism.
  
$$\xymatrix{K\ar@{^{(}->}[r]^k \ar@{->>}[d]^p&R\ar[r]_{r_0}\ar[d]^i&X\ar@/_1pc/[l]_{\delta}\ar@{=}[d]\\%
L\ar@{^{(}->}[r]^{\phi}&P\ar[r]_{p_0}&X\ar@/_1pc/[l]_{\sigma}}$$
\end{proof}

We now are able to prove our theorem :
\begin{theorem}
A category $\mathbb C$ is  finitely cocomplete homological if and only if it satisfies the three axioms A1-A3, and is semi-abelian   iff it satisfies all the four  axioms A1-A4. Moreover, in axiom A1, not all finite colimits are needed : finite sums and coequalizers of equivalence relations suffice.
\end{theorem}
\begin{proof}
Let A1' be the condition : 
 $\mathbb C$ is pointed and has finite products and coproducts, all equalizers, and coequalizers of equivalence relations.
We first prove that any category which satisfies A1', A2 and A3 is finitely cocomplete homological. Of course, condition A2 implies condition (F) in Proposition \ref{prodsemdir}, hence implies protomodularity, so in view of Proposition \ref{homolcoc} it suffices to prove that pullbacks of regular epimorphisms are regular epimorphisms ; hence in view of axiom A3, it  suffices to show that axioms A1' and A2 imply that any regular epimorphism is a cokernel. But by Proposition \ref{critproto}, $\mathbb C$ satisfies condition (A) of Lemma \ref{Manfred}, hence any regular epimorphism indeed is a cokernel.

Conversely, since by the same argument, in any pointed protomodular category any regular epimorphism is a cokernel, it follows that any finitely cocomplete homological category satisfies A1, A2 and A3.

Proposition \ref{homolex} then achieves the proof.
\end{proof}

\noindent{\bf Remark :} It should be noted that for a regular category, Axiom A4 is equivalent with:

\N{A4'.} In the following commutative diagram, if $f$ is a kernel, $g$  and $h$ are regular epimorphisms and $k$ is a monomorphism, then $k$ is a kernel.\medskip
$$\xymatrix{
A\ar[d]_g\ar[r]^f&B\ar[d]^h\\
C\ar[r]_k&D}$$
and of course if A1' and A2 hold, then ``coequalizer'' may be replaced by ``cokernel''.

As an easy consequence of our axioms we state the following result which does not seem to have been pointed out before :

\begin{theorem}\label{expo}If $I$ is a small category and if $\mathbb C$ is a pointed protomodular (resp homological, finitely cocomplete homological,  semi-abelian) category , then so is the category $\mathbb C^I$ of functors from $I$ to $\mathbb C$.\end{theorem}
\begin{proof}It is well-known that if limits of colimits of a specified type exist in $\mathbb C$, then they exist in $\mathbb C^I$, and are computed pointwise. Moreover, since finite limits exist  in $\mathbb C$ hence in $\mathbb C^I$, being a regular epimorphism can be rephrased as ``being the coequalizer of its kernel pair'', hence it is true for a natural transformation $\alpha$ in $\mathbb C^I$ if and only if it is pointwise true for each $\alpha_C$ in $\mathbb C$. Also, ``being a monomorphism'' can be verified pointwise, because it can be expressed in terms of finite limits (in any category, a morphism $f : A\rightarrow B$ is a monomorphism if and only if the pair $(1_A,1_A)$ is its kernel pair). And if cokernels exist, being a kernel may be rephrased as being the kernel of its cokernel. So if one  of the axioms A1, A2, A3 holds in $\mathbb C$, then it holds in $\mathbb C^I$. This is also true for axiom A4 (or more precisely its modified form A4') if a zero object, finite limits and finite colimits exist. So we have the result concerning finitely cocomplete homological and semi-abelian categories. For protomodular, homological and exact homological categories, the only axiom that needs some (very little !) work is protomodularity. But if $k: K\mono X$ is the kernel of a morphism $p: X\rightarrow Y$ with section $s:Y\rightarrow X$ in $\mathbb C^I$, and if $k$ and $s$ factor through some monomorphism $m:W\mono X$ by, say, $k=mk'$ and $s=ms'$, then this remains true pointwise, i.e. that each $k_C$ is the kernel in $\mathbb C$ of $p_C$, which has section $s_c$, and $k_c$ and $s_c$ factor through the monomorphism $m_C$ by $k'_c$ and $s'_c$ resp. So by protomodularity in $\mathbb C$ each $m_c$ is an isomorphism, hence so is $m$ in $\mathbb C^I$.
\end{proof}

 \small{\centerline{\sc\bf Acknowledgement:}\medskip
 
We express our gratitude  to Dominique Bourn for helpful discussions.

\end{document}